\documentclass[11pt]{amsart}
\usepackage{amssymb,amsfonts,amscd,amsbsy}
\usepackage[mathscr]{eucal}
\usepackage{url}

\theoremstyle{plain}

\theoremstyle{definition}

\begin{document}
\title{Principal forms $X^2+nY^2$ representing many integers}
\author{David Brink, Pieter Moree and Robert Osburn}

\address{School of Mathematical Sciences, University College Dublin, Belfield, Dublin 4, Ireland}

\address{Max-Planck-Institut f\"ur Mathematik, Vivatsgasse 7, D-53111 Bonn, Germany}

\email{david.brink@ucd.ie}

\email{moree@mpim-bonn.mpg.de}

\email{robert.osburn@ucd.ie}


\date{\today}

\begin{abstract}
In 1966, Shanks and Schmid investigated the asymptotic behavior of the number of positive integers less than or equal to $x$ which are represented by the quadratic form $X^2+nY^2$. Based on some numerical computations, they observed that the constant occurring in the main term appears to be the largest for  $n=2$. In this paper, we prove that in fact this constant is unbounded as $n$ runs through positive integers with a fixed number of prime divisors.
\end{abstract}

\maketitle
\section{Introduction}

It is a classical result of Landau \cite{L} from 1908 that the number $B(x)$ of integers less than or equal to $x$ which are representable as the sum of two squares $X^2+Y^2$ satisfies the asymptotic formula
\begin{equation} \label{lan}
\displaystyle B(x) \sim C \frac{x}{\sqrt{\log x}}\ \ \textrm{as}\ x \to \infty
\end{equation}
with the constant
\begin{equation*}
\displaystyle C=\frac{1}{\sqrt{2}} \prod_{p\equiv 3\ (\textrm{mod}\ 4)} \Biggl({1 \over {1 - 1/p^2}} \Biggr)^{1/2}=0.764223654,
\end{equation*}
where $p$ denotes a prime.
Independently, Ramanujan in his first letter to Hardy in 1913 stated essentially that
\begin{equation*}
\displaystyle B(x)\sim C \int_{2}^{x} \frac{dt}{\sqrt{\log t}}\ \ \textrm{as}\ x \to \infty.
\end{equation*}
Later claims by Hardy \cite{Hardy} that Ramanujan's integral did not give a better estimate of $B(x)$ than Landau's simpler formula were shown to be false by Shanks \cite{shanks}. The constant $C$ is now called the \emph{Landau-Ramanujan constant}.

Consider a primitive quadratic form $f(X,Y)=aX^2 + bXY + cY^2$ over $\mathbb{Z}$ with non-square discriminant $D=b^2-4ac$, and suppose $f$ is positive in case it is definite.
Let $B_{f}(x)$ be the number of positive integers less than or equal to $x$ which are representable by $f$.
Paul Bernays, a doctoral student of Landau's at G{\"o}ttingen, proved the following generalization of (\ref{lan}) in his 1912 thesis \cite[pages 59 and 115--116]{B}:
\begin{equation} \label{paulc}
\displaystyle B_{f}(x) \sim C(D) \frac{x}{\sqrt{\log x}}\ \ \textrm{as}\ x \to \infty
\end{equation}
with a non-zero constant $C(D)$ depending only on $D$. Thus $C(-4)$ is the Landau-Ramanujan constant.
Bernays did not explicitly give $C(D)$ for any other value of $D$. The problem of computing these constants has subsequently attracted considerable attention.

The original method of Landau can be used to compute $C(D)$ when the class number $h(D)$ is 1 or, with some additional complications, not too big. In 1966, Shanks and Schmid \cite{ss} studied the forms $f(X, Y)=X^2 + nY^2$ and determined the corresponding constants $C(-4n)$ ($b_n$ in their notation) in this way for 30 values of $n$ in the range $-34\leq n\leq 256$ with class numbers $h(-4n)=1$, 2, 3, 4, and 8. In particular they find
\[
C(-8)=0.872887558.
\]
They then state (page 561) ``We note, in passing, that of all binary forms $u^2 + nv^2$, $u^2 + 2v^2$ is the most populous, since $b_{2}$ is the largest of these constants.'' It is not completely clear as to whether they meant that $C(-8)$ is the largest amongst the values computed or that the maximum value of $C(-4n)$ as $n$ ranges over all integers is assumed for $n=2$. In any case, this quote motivates the following question: Is $C(-8)$ the maximum value? The purpose of the present paper is to answer this question.
Specifically we prove:
\\[5mm]
\textbf{Theorem 1.1}
\emph{If $\Delta$ is a fixed negative fundamental discriminant, then $C(\Delta q)$ is unbounded as $q$ runs through the primes congruent to 1 modulo 4.
If $\Delta$ is a fixed positive fundamental discriminant or $1$, then $C(-\Delta q)$ is unbounded as $q$ runs through the primes congruent to 3 modulo~4.}
\\[5mm]
It follows for example that $C(-4q)$ is unbounded as $q$ runs through the primes congruent to 1 modulo 4. However, it is not easy finding a concrete such $q$ with $C(-4q)>C(-8)$.
We have been able to find only one such example, namely
\begin{equation}\label{q}
C(-4\cdot 13779962790518414129 )= 0.875986.
\end{equation}

In Section \ref{twee} we present a formula for $C(D)$ in the case $D$ is a fundamental discriminant and sketch how it
is derived. In Section \ref{drie} we discuss computational aspects of this formula. In Section \ref{vier} we
prove Theorem 1.1 by making some adjustments in a proof of Joshi \cite{Joshi}. In the final section we
formulate some further problems and questions related to $C(D)$.

\section{Making $C(D)$ explicit}
\label{twee}

Using results of Kaplan and Williams \cite{kw} and Sun and Williams \cite{sw}, an explicit formula for Bernays' constant $C(D)$ for discriminants $D<0$ was given in \cite{mo} (see (2.5), (2.8), and (2.11)).
When $D$ is a fundamental discriminant, i.e., equal to the discriminant of a quadratic number field, this formula reduces to
\begin{equation}
\label{bern}
\displaystyle C(D)=\frac{1}{2^{\omega(D)-1}} \Biggl [ \frac{|D|}{\varphi(|D|)} \frac{L(1, \chi_D)}{\pi} E(D) \Biggr ]^{1/2}.
\end{equation}
Here $\omega(D)$ is the number of prime divisors in $D$, $\varphi$ is Euler's phi function, $L(\ \cdot\ ,\chi_D)$ is the Dirichlet L-series corresponding to the Kronecker symbol $\chi_D=(D/\ \cdot\ )$, and $E(D)$ is the Euler product
\[
E(D)=\prod_{(D/p)=-1}{1\over 1-1/p^2}.
\]
Note that
\[
1< E(D)<\prod_{p}\frac{1}{1-1/p^2}=\frac{\pi^2}{6}
\]
and thus the contribution of $E(D)$ to $C(D)$ is limited.
 
\indent We will now sketch some arguments that go into the derivation of (\ref{bern}). Recall that a primitive positive
definite quadratic form $aX^2+bXY+cY^2$, $[a,b,c]$ for short, is said to be {\it reduced} if
$|b|\le a\le c$, and $b\ge 0$ if either $|b|=a$ or $a=c$. Every primitive positive definite form is properly equivalent to a unique reduced form \cite[Theorem 2.8]{cox}. This reduced form represents precisely the same integers as the original one. Thus we might assume from the outset that our form is reduced. We say that two forms are in the same {\it class} if they are properly equivalent. It is easy to see that the number $h(D)$ of classes of primitive positive definite forms of discriminant $D$ is finite, and furthermore $h(D)$ is equal to
the number of reduced forms of discriminant $D$ \cite[Theorem 2.13]{cox}. In the case $D$ is a fundamental discriminant, $h(D)$ also equals the class number of the quadratic number field $\mathbb Q(\sqrt{D})$, and we have
\begin{equation*} \label{h}
{L(1,\chi_D)\over \pi}=u(D)\cdot {h(D)\over  \sqrt{-D}}
\end{equation*}
where $u(-3)=1/3$, $u(-4)=1/2$, and $u(D)=1$ for $D<-4$.\\
\indent We say that two primitive positive definite forms of discriminant $D$ are
in the same {\it genus} if they represent the same values in $(\mathbb Z/D\mathbb Z)^*$. The number of genera is known to be a power of two and so can be written as $2^{t(D)}$ for some $t(D)\ge 0$. The number $t(D)$ depends only on $\omega(D)$ and the the congruence class of $D$ modulo 32 (see (2.7) in \cite{mo}). Furthermore, given a negative integer $D \equiv 0$, $1 {\pmod 4}$, the {\it principal form} is defined by $\bigl[1, 0, -\frac{D}{4} \bigr]$ if $D \equiv 0 \pmod 4$ and 
$\bigl[1, 1, \frac{1-D}{4} \bigr]$ if $D \equiv 1 \pmod 4$. The principal form has discriminant $D$ and is reduced. If $D=-4n$, then $[1, 0, n]$ is the principal form. Finally, we say that an integer $m$ is represented by the genus $G$ if it is represented by at least one class in $G$. For example, if $n=14$ (and thus $D=-56$ and $h(D)=4$), then we have two genera $G_1$ and $G_2$ where, say, (the class of) $[1, 0, 14]$ and $[2, 0, 7]$ belong to $G_1$ while $[3, -2, 5]$ and $[3, 2, 5]$ belong to $G_2$. 

\indent Let us first consider the simpler problem of deriving an analogue of (\ref{paulc}) with
$B_f'(x)$ instead of $B_f(x)$, where $B_f'(x)$ counts the number of integers $m\le x$ represented by $f$ with $(m,D)=1$.
Note that without loss of generality we may assume that $f$ is reduced. Let $f$ belong to the genus $G$. The counting strategy is
as follows:\\
\begin{enumerate}
\item Compute $B_D'(x)$; the number of integers $m\le x$ coprime to $D$ that are represented by {\it any}
reduced form of
discriminant $D$.\\
\item Compute $B_G'(x)$; the number of integers $m\le x$ coprime to $D$ that are represented by $G$. \\
\item Compute $B_f'(x)$,
\end{enumerate}
where the result of each step provides input for the next. We now proceed through these steps.

\indent (1) An integer $m\le x$ is counted if and only if its prime divisors that occur to an odd power satisfy
$(D/p)=1$. The associated L-series $L_D(s)=\sum m^{-s}$ thus has, for Re$(s)>1$, the following Euler product:
\begin{equation*}
L_D(s)=\prod_{(D/p)=1}(1-p^{-s})^{-1}\prod_{(D/p)=-1}(1-p^{-2s})^{-1}. 
\end{equation*}
Using the Euler product formula for $\zeta(s)$ and $L(s,\chi_D)$ one finds, for Re$(s)>1$, that
\begin{equation*}
L_D(s)^2=\zeta(s)L(s,\chi_D)\prod_{p|D}(1-p^{-s})\prod_{(D/p)=-1}(1-p^{-2s})^{-1}=\zeta(s)g(s),
\end{equation*}
say. Then using the Selberg-Delange method (see, e.g., \cite[Proposition 5]{FMS}), one finds
that
$$B_D'(x)\sim \sqrt{g(1)\over \pi}{x\over \sqrt{\log x}},$$
where
$$\sqrt{g(1)\over \pi}=\Biggl [ \frac{\varphi(|D|)}{|D|} \frac{L(1, \chi_D)}{\pi} E(D) \Biggr ]^{1/2} =: J(D),$$
the James constant, who first established this result \cite{James} (from a modern perspective this result
is completely standard). For different proofs see \cite{mo, williams}. The $\sqrt{\log x}$ is due to the fact
that asymptotically half of the primes are represented by a reduced binary quadratic form. If we have many
small primes $p$ with $(D/p)=1$, then $L(1,\chi_D)$ is large on the one hand, and on the other many integers
$m$ are counted by $B_D'(x)$. Thus it is natural to expect that $J(D)$ scales as an increasing function
of $L(1,\chi_D)$.\\
(2) Any integer coprime to $D$ is represented by at most one genus. Bernays showed that, as $x$ tends to
infinty, the integers $m\le x$ become equidistributed over the genera, that is, we have
$$B_G'(x)\sim {B_D(x)\over 2^{t(D)}}\sim {J(D)\over 2^{t(D)}}{x\over \sqrt{\log x}}.$$
(3) We let $E$ be the set of integers that are represented by $G$, but not by $f$, and $E(x)$ the associated
counting function (for example, in the above example $2$ and $7$ are represented by $G_1$, but not by $[1,0,14]$).
Bernays showed that $E(x)=o(x\log^{-1/2}x)$
(this result was later sharpened by Fomenko \cite{Fomenko} to $E(x)\ll x\log^{-2/3}x$). We finally conclude that
$$B_f'(x)\sim B_G'(x)\sim {J(D)\over 2^{t(D)}}{x\over \sqrt{\log x}}.$$
This solves the asymptotic counting problem for $B_f'(x)$. However, we are interested in $B_f(x)$. A complication that arises
here is that more than one genus might represent an integer not coprime to $D$. Bernays took this complication into account and arrived at the following formula
\begin{equation*}
\label{bernabas}
C(D)={J(D)\over 2^{t(D)}}\sum_{m|D^{\infty}}{g(m,D)\over m},
\end{equation*}
where $g(m,D)$ denotes the number of genera  of discriminant $D$ representing $m$ and
$m|D^{\infty}$ means that $m$ divides some arbitrary power of $D$. The latter sum is a rational
number and its explicit evaluation was only made possible by the recent papers
\cite{kw, sw} mentioned in the beginning of this section. In the case $D$ is a fundamental discriminant
the latter sum equals $|D|/\varphi(|D|)$, $t(D)=\omega(D)-1$ and we obtain (\ref{bern}). For the general
formula we refer to \cite{mo}.\\

\noindent {\tt Remark 1}. In retrospect one sees why the classical Landau case readily follows. There
one has $B_{-4}'(x)=B_G'(x)=B_f'(x)$.
Let $Od(x)$ denote the number of odd integers $\le x$ that can be written as a sum of two squares.
Thus $B_{-4}'(x)=Od(x)$. Note that $J(-4)=C/2$ and $B(x)=\sum_{j=0}^{\infty}Od(x2^{-j})$. We infer
that $B(x)\sim 2Od(x)\sim Cx \log^{-1/2} x$. Observe that $h(-4)=1$.
Indeed, the three cases coincide if and only if $h(D)=1$. The second and third case coincide (that is every genus consists of one class)
if and only if $n\ge 1$ is a convenient number (``numerus idoneus"). This happens for at least 65 and at
most 66 integers \cite[pp. 61-62]{cox}.

\section{Computation of $C(D)$}
\label{drie}
The class numbers $h(D)$ can be computed for about $|D|<10^{20}$ by GP/PARI. In order to compute $C(D)$ accurately, the problem thus lies in computing the Euler product $E(D)$. Taking the product over, say, the first 100,000 primes gives a precision of about 6 decimals. A much better precision is obtained by the rapidly converging infinite product (see \cite[eq.\ (16)]{shanks} or \cite[eq.\ (3.3)]{mo})
\begin{equation*}\label{EP}
E(D)=\prod_{k=1}^\infty\left({\zeta(2^k)\over L(2^k,\chi_D)}\prod_{p|D}\left(1-p^{-2^k}\right)\right)^{1/2^k}.
\end{equation*}
However, the values $L(2^k,\chi_D)$ can be computed only for about $|D|<10^6$ by GP/PARI .

Using the above formulas, we compute $C(D)$ for some small fundamental discriminants $D<0$ ordered according as to whether $\omega(D)$ is 1, 2, 3, or 4:
\[
\begin{array}[t]{r|l}
D & \ \ \ \ C(D)\\
\hline
-3 & 0.638909405 \\
-4 & 0.764223654 \\
-7 & 0.724719521 \\
-8 & 0.872887558 \\
-11 & 0.677388018 \\
-19 & 0.606300131 \\
-23 & 0.841512352 \\
-31 & 0.801014576 \\
-43 & 0.500610055 \\
-47 & 0.891550880 \\
-59 & 0.735485997 \\
-67 & 0.448813095 \\
-71 & 0.938541302 \\
-79 & 0.812629337 \\
-83 & 0.684502354 \\
\end{array}
\ \
\begin{array}[t]{r|l}
D & \ \ \ \ C(D)\\
\hline
-15 & 0.501918636 \\
-20 & 0.535179999 \\
-24 & 0.558357114 \\
-35 & 0.407379938 \\
-39 & 0.518747305 \\
-40 & 0.473558100 \\
-51 & 0.390646647 \\
-52 & 0.420720518 \\
-55 & 0.458949554 \\
-56 & 0.563486772 \\
-68 & 0.520288297 \\
-87 & 0.512573818 \\
-88 & 0.375792661 \\
-91 & 0.317487516 \\
-95 & 0.528624390 \\
\end{array}
\ \
\begin{array}[t]{r|l}
D & \ \ \ \ C(D)\\
\hline
-84 & 0.310647641 \\
-120 & 0.296417662 \\
-132 & 0.274765289 \\
-168 & 0.267006498 \\
-195 & 0.220993565 \\
-228 & 0.237562625 \\
-231 & 0.309699577 \\
-255 & 0.307681243 \\
-260 & 0.293752522 \\
-264 & 0.319941656 \\
-276 & 0.309309571 \\
-280 & 0.223644570 \\
-308 & 0.277034255 \\
-312 & 0.223049066 \\
-340 & 0.204812008 \\
\end{array}
\ \
\begin{array}[t]{r|l}
D & \ \ \ \ C(D)\\
\hline
-420 & 0.164080141 \\
-660 & 0.143806822 \\
-840 & 0.139069358 \\
-1092 & 0.123274604 \\
-1140 & 0.171607125 \\
-1155 & 0.109195133 \\
-1320 & 0.121504603 \\
-1380 & 0.117420083 \\
-1428 & 0.114424422 \\
-1540 & 0.108139197 \\
-1560 & 0.161366493 \\
-1716 & 0.148895032 \\
-1848 & 0.109066658 \\
-1860 & 0.151207258 \\
-1995 & 0.093833104 \\
\end{array}
\]
It appears that $2^{1-\omega(D)}$ tends to dominate the other factors in (\ref{bern}).
We have $\omega(D)=1$ when $D$ is of the form $-q$ for a prime $q\equiv 3\ (\textrm{mod}\ 4)$. It is straightforward to find such primes with $C(-q)>C(-8)$, for example $q=47$, $71$, $167$, $191$ or $239$. The largest value of $C(-q)$ with $q<10^9$, which is also the largest value of $C(D)$ that we know, is
\[
C(-984452999) = 1.527855.
\]
But already for $\omega(D)=2$ it becomes much more difficult to find a $D$ with $C(D)>C(-8)$, and the only example we know is (\ref{q}).
The difficulty in finding such examples is explained by the fact that
\[
\max_{0<-D<x}L(1,\chi_D)
\]
grows very slowly with $x$. For example, Bateman, Erd{\"o}s, and Chowla \cite{BCE} proved
$$
L(1, \chi_D) < \frac{10}{3} \frac{\varphi(|D|)}{|D|} \log |D| + 1.
$$
Moreover, assuming a suitable generalized Riemann hypothesis, Littlewood \cite{little} showed
\begin{equation} \label{Litt}
e^\gamma\leq \limsup_{D\to-\infty} \frac{L(1,\chi_D)}{\log\log |D|}\leq 2e^\gamma
\end{equation}
with $D$ running through negative fundamental discriminants. The left-hand inequality in (\ref{Litt}) was shown unconditionally by Chowla \cite{Chowla} (see also the discussion in \cite{BCE}). Recent work by Granville and Soundararajan \cite{Granville} gives strong evidence via a probabilistic model that $e^\gamma$ is in fact the true limit superior of $L(1,\chi_D)/\log\log |D|$.

Regarding small values of $C(D)$, we mention that it was shown in \cite{mo} that $C(D)\cdot\sqrt{-D}$ is minimal for $D=-3$ as $D$ ranges over all negative discriminants, and that accordingly of all the two-dimensional lattices of covolume $1$, the hexagonal lattice has asymptotically the fewest distances.

\section{Proof of Theorem 1.1}
\label{vier}
In order to prove Theorem 1.1, we need to show that $L(1, \chi_D)$ is unbounded when $D$ runs through a certain subset of discriminants. In this direction, it is proved in \cite{BCE} that
\begin{equation} \label{ok}
\limsup_{D\to-\infty} \frac{L(1,\chi_D)}{\log\log |D|}\geq \frac{e^\gamma}{18}
\end{equation}
where $D$ runs through fundamental discriminants of the form $D=-q$ with $q\equiv 3\ (\textrm{mod}\ 4)$ prime. This implies that $C(D)$ is unbounded, but says nothing about discriminants of the form $D=-4n$. Our main interest is in a result of Joshi \cite{Joshi} in which she improved (\ref{ok}) by removing the factor 18 (for a quantitative version of this result, see \cite{cook}). It turns out that one can make suitable adjustments to Joshi's proof in order to prove the following:
\\[5mm]
\textbf{Theorem 3.1.}
\emph{Let $\Delta$ be a fundamental discriminant or 1, let $c$ and $d$ be coprime integers with $d$ divisible by $\Delta$ and 8,  let $q$ run through the primes congruent to $c\ (\textrm{mod}\ d)$, and let $\chi$ be the Kronecker character $(\Delta q^*/\  \cdot\ )$ with $q^*=\lambda q$, $\lambda=(-1)^{(c-1)/2}$. Then
\[
\limsup_{\substack{q\to\infty\\q\equiv c\ (\textrm{mod}\ d)}}\frac{L(1,\chi)}{\log\log q}
\geq e^\gamma\cdot\prod_{p|d}
\frac{1-\frac{1}{p}}
{1-\left(\frac{\Delta c^*}{p}\right)\frac{1}{p}},
\]
where $c^*=\lambda c$ and $\gamma$ is Euler's constant.}
\\[5mm]
Note that this is close to being best possible, cf.\ (\ref{Litt}).
\\[5mm]
\emph{Proof of Theorem 3.1.}
The theorem is a generalization of \cite[Theorem 2]{Joshi} which corresponds to the case $\Delta=1$. As the proof is a modification of Joshi's argument, we give the necessary changes.

Fix some (small) $\epsilon>0$. It suffices to show that for every (large) $x$ there exists a prime $q\leq x$, $q\equiv c\ (\textrm{mod}\ d)$, such that
\begin{eqnarray} \label{start}
\nonumber
\log L(1,\chi)&\geq& \log\log\log x+\gamma+\sum_{p|d}\log\left(1-\frac{1}{p}\right)\\
&&-\sum_{p|d}\log\left(1-\left(\frac{\Delta c^*}{p}\right)\frac{1}{p}\right)+\log(1-2\epsilon)+o(1).
\end{eqnarray}

We prove (\ref{start}) by constructing a set $\Sigma=\Sigma(x)$ of primes $q\leq x$, $q\equiv c\ (\textrm{mod}\ d)$, with $S=|\Sigma|$ and showing
\begin{eqnarray} \label{key}
\nonumber
\sum_{q\in\Sigma} \log L(1,\chi)&\geq &
S\left(\log\log\log x+\gamma+\sum_{p|d}\log\left(1-\frac{1}{p}\right)\right.\\
&&\left.-\sum_{p|d}\log\left(1-\left(\frac{\Delta c^*}{p}\right)\frac{1}{p}\right)
+\log(1-2\epsilon)\right)+o(S).
\end{eqnarray}

Put
\[
y=(\log x)^{1-2\epsilon}
\]
and let $p_1,\dots,p_m$ be the primes not greater than $y$ and not dividing $d$.
Define $r$ as in \cite[p.\ 64]{Joshi}, and let
$k=dp_1\cdots p_{r-1}p_{r+1}\cdots p_m$. For each $i\neq r$, let $g_i$ (respectively $h_i$) be a quadratic residue (respectively non-residue) modulo $p_i$. Let $l\leq k$ be the unique positive integer satisfying $l \equiv c\ (\textrm{mod\ } d)$ and
\[
l \equiv \begin{cases}
g_i\ (\textrm{mod\ }p_i)& \text{for $\bigl({\lambda \Delta \over p_i}\bigr)=1$, $i\neq r$} \\
h_i\ (\textrm{mod\ }p_i)& \text{for $\bigl({\lambda \Delta \over p_i}\bigr)=-1$, $i\neq r$}.
\end{cases}
\]
Define
\[
\Sigma=\{q\ \textrm{prime}\ |\ \sqrt{x}\leq q\leq x,\ q\equiv l\ (\textrm{mod\ }k)\}.
\]
Then every $q\in\Sigma$ satisfies $q\equiv c\ (\textrm{mod\ }d)$ and $\chi(p_i)=1$ for  $i\neq r$ since
\[
\chi(p_i)=
\left(\dfrac{\Delta q^*}{p_i}\right)=
\left(\dfrac{\lambda\Delta}{p_i}\right)
\left(\dfrac{q}{p_i}\right)=
\left(\dfrac{\lambda\Delta}{p_i}\right)
\left(\dfrac{l}{p_i}\right)=1.
\]

So far, the only difference compared with Joshi's proof is the definition of $l$ and $\chi$ (in \cite{Joshi}, $\chi$ is the character $\bigl({\cdot \over q}\bigr)=\bigl({q^* \over \cdot}\bigr)$ corresponding to $\Delta=1$). The different definition of $l$  plays no role other than guaranteeing that we still have $\chi(p_i)=1$, cf.\ \cite[p.\ 65]{Joshi}.
Hence, as in \cite[(24)]{Joshi}, we get
\begin{eqnarray*}
\nonumber
\sum_{q\in\Sigma} \log L(1,\chi)&\geq &
S\left(\log\log\log x+\gamma+\sum_{p|d}\log\left(1-\frac{1}{p}\right)\right.\\
\label{2}
&&\left.-\sum_{p|d}\log\left(1-\left(\frac{\Delta c^*}{p}\right)\frac{1}{p}\right)
+\log(1-2\epsilon)\right)+R+o(S)
\end{eqnarray*}
where
\begin{equation*}
R=\sum_{q\in\Sigma}\sum_{p>y}\frac{\chi(p)}{p}.
\end{equation*}
We now show $R=o(S)$ and hence (\ref{key}) by splitting the summation over $p$ into five intervals $I_1,\dots,I_5$ and thus writing $R=R_1+\dots+R_5$ with
\begin{equation*}
R_i=\sum_{q\in\Sigma}\sum_{p\in I_i}\frac{\chi(p)}{p}.
\end{equation*}

The estimation of $R_1$ and $R_2$ is practically the same as in
Joshi's paper, only one has to replace $\bigl({\lambda \over p}\bigr)$ by  $\bigl({\lambda \Delta \over p}\bigr)$ in \cite[(27)]{Joshi} and the equation below that, which makes no difference since the sign of that factor plays no role anyway. The estimation of $R_3$ is exactly the same since it relies on the majorization
\[
\left|\sum_{p\in I_3}\frac{\chi(p)}{p}\right|\leq \sum_{p\in I_3}\frac{1}{p}.
\]

The estimation of $R_4$ requires some more care since it relies on the large sieve as stated in \cite[Lemma 1]{Joshi} which works only for prime moduli. Put $\beta=2+\epsilon^{-1}$ and subdivide $I_4$ into intervals $J_t$ each containing $Z_t$ primes as in \cite[p.\ 70]{Joshi}. Then \cite[(30)]{Joshi} remains valid, i.e.
\begin{equation}
\label{O}
\sum_{q\in \Sigma}\sum_{p\in J_t}\frac{\chi(p)}{p}
-\frac{1}{t}\sum_{q\in \Sigma}\sum_{p\in J_t}\chi(p)
=O\left(\frac{S}{(\log x)^{2\beta}}\right).
\end{equation}
Let $J_t^+$ and $J_t^-$ denote the sets of primes in $J_t$ with $\bigl({\Delta \over p}\bigr)=1$ and $\bigl({\Delta \over p}\bigr)=-1$, respectively.
Then $Z_t=Z_t^+ + Z_t^-$ where
$Z_t^+$ and $Z_t^-$ are defined analogously.
Also, let $Z_t(a,q)$ be the number of $p$ in $J_t$ which are congruent to $a$ modulo $q$, and similarly write $Z_t(a,q)=Z_t^+(a,q)+Z_t^-(a,q)$.
Then a computation using the large sieve, cf. \cite[p. 71]{Joshi}, shows
\begin{eqnarray*}
\nonumber
\left|\frac{1}{t}\sum_{q\in \Sigma}\sum_{p\in J_t^+}
\left(\dfrac{p}{q}\right)
\right|^2
&=&
\frac{1}{t^2}\left|\sum_{q\in \Sigma}\sum_{j=1}^{q-1}
\left(\dfrac{j}{q}\right)
\left(Z_t^+(j,q)-\frac{Z_t^+}{q}\right)\right|^2\\
&\leq&
\frac{S^2}{(\log x)^{4\beta}},
\end{eqnarray*}
and similarly with the summation over $p\in J_t^-$. Since $\chi(p)=\bigl({\Delta \over p}\bigr)\bigl({p \over q}\bigr)$,
we now get
\begin{eqnarray}
\nonumber
\left|\frac{1}{t}\sum_{q\in \Sigma}\sum_{p\in J_t}\chi(p)\right|
&\leq&
\left|\frac{1}{t}\sum_{q\in \Sigma}\sum_{p\in J_t^+}\chi(p)\right|+
\left|\frac{1}{t}\sum_{q\in \Sigma}\sum_{p\in J_t^-}\chi(p)\right|\\
\nonumber
&=&
\left|\frac{1}{t}\sum_{q\in \Sigma}\sum_{p\in J_t^+}\left(\dfrac{p}{q}\right)\right|+
\left|\frac{1}{t}\sum_{q\in \Sigma}\sum_{p\in J_t^-}\left(\dfrac{p}{q}\right)\right|\\
\nonumber
&\leq&
2\cdot \frac{S}{(\log x)^{2\beta}}\\
\label{OO}
&=&
O\left(\frac{S}{(\log x)^{2\beta}}\right).
\end{eqnarray}
{}From \eqref{O} and \eqref{OO} follows
\[
\left|\sum_{q\in \Sigma}\sum_{p\in J_t}\frac{\chi(p)}{p}\right|
=O\left(\frac{S}{(\log x)^{2\beta}}\right),
\]
and thus
\[
R_4=\sum_{q\in \Sigma}\sum_{p\in I_4}\frac{\chi(p)}{p}=o(S).
\]
Finally, the estimation of $R_5$ can be carried out as in \cite[p. 72]{Joshi} by writing
\[
\sum_{v<p\leq w}\frac{\chi(p)}{p}=
\sum_{j=1}^{|\Delta q|} \chi(j) 
\sum_{\substack{v<p\leq w\\p\equiv j\ (\textrm{mod\ } |\Delta q|)}}\frac{1}{p}
\]
and using \cite[Lemma 3]{Joshi}. \qed
\\[5mm]
{\it Proof of Theorem 1.1.} If $\Delta$ is a negative fundamental discriminant, then let $c=1$ and $d$ be divisible by $\Delta$ and $8$. Hence Theorem 3.1 implies that for $\chi=\chi_{D}$, $D=\Delta q$, we have
$$
\sup_{\substack{q\to\infty\\ q\equiv c\ (\textrm{mod\ }d)}} L(1, \chi)=\infty.
$$
Applying this to (\ref{bern}) yields the first statement. If $\Delta$ is a positive fundamental discriminant or $1$, let $c=-1$ and $d$ be divisible by $\Delta$ and $8$. Applying Theorem 3.1 with $D=-\Delta q$ to (\ref{bern}) implies the second statement. \qed

\section{Outlook}
We focused on large values of $C(\Delta q)$. Likewise, one might ask about small values. Is it true
for example that $\lim \inf C(\Delta q)=0$ as $q$ runs over the primes congruent to 1 modulo 4?
Is $C(\Delta q)$ constant on average? Likewise one can consider our table and wonder whether the values
in a column are constant on average. If so, will this constant be zero or not? Perhaps variations of the
techniques in \cite{Granville} can be used to study this. A further open problem is to determine whether $B_f(x)$ is asymptotically better approximated
by $C(D)\int_2^x{dt\log^{-1/2}t}$ or $C(D)x\log^{-1/2}x$. Finally, it might be of interest to recover $C(D)$ for any discriminant $D<0$ using sieve methods. For example, the Landau-Ramanujan constant $C$ was verified (see \cite{iwaniec}) using the half-dimensional sieve.

\section*{Acknowlegements}
The first author was supported by grant 272-08-0323 from the Danish Agency for Science, Technology and Innovation while the third author was partially supported by Science Foundation Ireland 08/RFP MTH1081. The third author would like to thank the Institut des Hautes {\'E}tudes
Scientifiques for their hospitality and support during the preparation of this paper and Greg Martin and
Olivier Ramar\'e for their comments.

\end{document}